\documentstyle[10pt]{amsart}
\numberwithin{equation}{section}
\input{amssym.def}
\input{amssym.tex}

\begin{document}

\title[Mazur-Ulam theorem in non-Archimedean $n$-normed spaces]
 {A Mazur-Ulam theorem for Mappings of conservative distance in non-Archimedean $n$-normed spaces}

\author[Hahng-Yun Chu and Se-Hyun Ku]{Hahng-Yun Chu$^\dagger$ and Se-Hyun Ku$^\ast$}

\address{Hahng-Yun Chu, Department of Mathematical Sciences, Korea Advanced Institute of Science and Technology,
335, Gwahak-ro, Yuseong-gu, Daejeon 305-701, Republic of Korea }
\email{\rm hychu@@kaist.ac.kr}

\address{Se-Hyun Ku, Department of Mathematics, Chungnam National University, 79, Daehangno, Yuseong-Gu, Daejeon 305-764, Republic of Korea} \email{\rm shku@@cnu.ac.kr}

\thanks{
\newline\ \noindent 2010 {\it Mathematics Subject Classification}. primary 46B20, secondary 51M25, 46S10
\newline {\it Key words and Phrases}. Mazur-Ulam theorem, $n$-isometry, non-Archimedean $n$-normed space
\newline{\it $\ast$ Corresponding author}
\newline{\it $\dagger$ The author was supported by the second stage of the Brain Korea 21 Project, The Development Project of Human Resources in Mathematics, KAIST in 2008. \\}}

\maketitle

\begin{abstract}
In this article, we study the notions of $n$-isometries in non-Archimedean $n$-normed spaces over linear ordered non-Archimedean fields, and prove the Mazur-Ulam theorem in the spaces.
Furthermore, we obtain some properties for $n$-isometries in non-Archimedean $n$-normed spaces.
\end{abstract}

\baselineskip=18pt

%=============================================================================
%      define words
%=============================================================================
%\newcommand{\nil}{{\rm nilrad}}
%\newcommand{\deg}{{\rm deg}}
%\renewcommand{\theenumii}{\alph{enumii}}

%=============================================================================
%      define theorems
%=============================================================================

\theoremstyle{definition}
  \newtheorem{df}{Definition}[section]
    \newtheorem{rk}[df]{Remark}
\newtheorem{ma}[df]{Main Theorem}
  \newtheorem{cj}[df]{Conjecture}
  \newtheorem{pb}[df]{Problem}
\theoremstyle{plain}
  \newtheorem{lm}[df]{Lemma}
  \newtheorem{eq}[df]{equation}
  \newtheorem{thm}[df]{Theorem}
  \newtheorem{cor}[df]{Corollary}
  \newtheorem{pp}[df]{Proposition}

\setcounter{section}{0}
%\begin{equation}\end{equation}\begin{lemma}\end{lemma}\begin{theorem}\end{theorem}
%\begin{eqnarray*}\end{eqnarray*}

%====================================================================================
\section{Introduction}
%====================================================================================

Let $X$ and $Y$ be metric spaces with metrics $d_X$ and $d_Y,$ respectively.
A map $f:X\rightarrow Y$ is called an isometry if $d_Y(f(x),f(y))=d_X(x,y)$ for every $x,y\in X.$
Mazur and Ulam\cite{mu32} were the first to treated the theory of isometry.
They have proved the following theorem;

\textbf{Mazur-Ulam Theorem}\,\,\,\textit{Let $f$ be an isometric transformation from a real normed vector space
$X$  onto a real normed vector space $Y$ with $f(0)=0$. Then $f$ is linear.}

It was a natural ask if the result holds without the onto assumption.
Asked about this natural question, Baker\cite{b71} answered that every isometry of a real normed linear space
into a strictly convex real normed linear space is affine.
The Mazur-Ulam theorem has been widely studied by \cite{j01,ms08,r01,rs93,rw03,v03,x01}.

Chu et al.\cite{cpp04} have defined the notion of a $2$-isometry which is suitable to represent the concept of an area preserving mapping in linear $2$-normed spaces.
In \cite{c07}, Chu proved that the Mazur-Ulam theorem holds in linear $2$-normed spaces
under the condition that a $2$-isometry preserves collinearity.
Chu et al.\cite{ckk08} discussed characteristics of $2$-isometries.
In \cite{as}, Amyari and Sadeghi proved the Mazur-Ulam theorem in non-Archimedean $2$-normed spaces under the condition of strictly convexity.
Recently, Choy et al.\cite{cck} proved the theorem on non-Archimedean $2$-normed spaces over linear ordered non-Archimedean fields without the strictly convexity assumption.

Misiak\cite{m89-1, m89-2} defined the concept of an $n$-normed space and investigated the space.
Chu et al.\cite{clp04}, in linear $n$-normed spaces, defined the concept of an $n$-isometry that is suitable to represent the notion of a volume preserving mapping.
In~\cite{cck09}, Chu et al. generalized the Mazur-Ulam theorem to $n$-normed spaces.
In recent years, Chen and Song\cite{cs} characterized $n$-isometries in linear $n$-normed spaces.

In this paper, without the condition of the strictly convexity, we prove the (additive) Mazur-Ulam theorem on non-Archimedean $n$-normed spaces.
Firstly, we assert that an $n$-isometry $f$ from a non-Archimedean space to a non-Archimedean space 
preserves the midpoint of a segment under some condition about the set of all elements of a valued field 
whose valuations are 1.
Using the above result, we show that the Mazur-Ulam theorem on non-Archimedean $n$-normed spaces over linear ordered non-Archimedean fields.
In addition, we prove that the barycenter of a triangle in the non-Archimedean $n$-normed spaces is $f$-invariant under different conditions from those referred in previous statements. 
And then we also prove the (second typed) Mazur-Ulam theorem in non-Archimedean $n$-normed spaces under some different conditions.

\medskip

%%%%%%%%%%%%%%%%%%%%%%%%%%%%%%%%%%%%%%%%%%%%%%%%%%%%%%%%%%%%%%%%%%%%%%%%%%%%%%%%%%%%%%%%%%%%%%%%%%%%%%%%
%%%%%%%%%%%%%%%%%%%%%%%%%%%%%%%%%%%%%%%%%%%%%%%%%%%%%%%%%%%%%%%%%%%%%%%%%%%%%%%%%%%%%%%%%%%%%%%%%%%%%%%%
%%%%%%%%%%%%%%%%%%%%%%%%%%%%%%%%%%%%%%%%%%%%%%%%%%%%%%%%%%%%%%%%%%%%%%%%%%%%%%%%%%%%%%%%%%%%%%%%%%%%%%%%
\section{The Mazur-Ulam theorem on non-Archimedean $n$-normed spaces}

In this section, we introduce a non-Archimedean $n$-normed space which is a kind of a generalization of a non-Archimedean $2$-normed space and we show the (additive) Mazur-Ulam theorem for an $n$-isometry $f$ defined on a non-Archimedean $n$-normed space, that is, $f(x)-f(0)$ is additive.
Firstly, we consider some definitions and lemmas which are needed to prove the theorem.

%=======================================================================================================

Recall that a {\emph{non-Archimedean}} (or {\emph{ultrametric}}) {\emph{valuation}} is given by a map $|\cdot|$ from a field ${\mathcal{K}}$ into $[0,\infty)$ such that for all $r,s\in{\mathcal{K}}$:

(i) $|r|=0$ if and only if $r=0$;

(ii) $|rs|=|r||s|$;

(iii) $|r+s|\leq \max\{|r|, |s|\}.$

\medskip

\medskip
If every element of ${\mathcal{K}}$ carries a valuation then a field ${\mathcal{K}}$ is called a {\emph{valued field}}, for a convenience,
simply call it a field. It is obvious that $|1|=|-1|=1$ and $|n| \leq 1$ for
all $n \in {\mathbb N}$. A trivial example of a non-Archimedean valuation
is the map $|\cdot|$ taking everything but $0$ into $1$ and
$|0|=0$ (see~\cite{nb81}).

\medskip

\medskip

A {\emph{non-Archimedean norm}} is a function $\| \cdot \| :{\mathcal X} \to [0, \infty)$ such that for all $r \in {\mathcal K}$ and $x,y \in {\mathcal X}$:

(i) $\|x\| = 0$ if and only if $x=0$;

(ii) $\|rx\| = |r| \|x\|$;

(iii) the strong triangle inequality
$$\| x+ y\| \leq \max \{\|x\|, \|y\|\}.$$
Then we say $({\mathcal X}, \|\cdot\|)$ is a {\emph{non-Archimedean space}}.

%======================================= definition ===============================================
\begin{df}\label{df31}
Let ${\mathcal X}$ be a vector space with the dimension greater than $n-1$ over a valued field $\mathcal{K}$
with a non-Archimedean valuation $|\cdot|.$
A function $\| \cdot, \cdots, \cdot \|:{\mathcal{X}}\times\cdots\times{\mathcal{X}}\rightarrow[0,\infty)$
is said to be a {\emph{non-Archimedean $n$-norm}} if

$(\textrm{i}) \ \ \|x_1, \cdots, x_n \|=0 \Leftrightarrow x_1,
\cdots, x_n
 \textrm{ are linearly dependent} ;$

$(\textrm{ii}) \ \ \|x_1, \cdots, x_n \| = \| x_{j_1}, \cdots,
x_{j_n} \| $ for every permutation $(j_1, \cdots, j_n)$ of $(1,
\cdots, n) ;$

$(\textrm{iii}) \ \ \| \alpha x_1, \cdots, x_n \| =| \alpha | \
\| x_1, \cdots, x_n \| ;$

$(\textrm{iv}) \ \ \|x+y, x_2, \cdots, x_n \| \le \max\{\|x, x_2,
\cdots, x_n\|, \|y, x_2, \cdots, x_n\|\} $

\noindent for all $\alpha \in \mathcal{K}$ and all $x, y, x_1,
\cdots, x_n \in \mathcal{X}$.
Then $(\mathcal{X},\| \cdot, \cdots, \cdot \|)$ is called a {\it non-Archimedean
$n$-normed space}.

\end{df}

From now on, let $\mathcal{X}$ and $\mathcal{Y}$ be non-Archimedean $n$-normed spaces
over a linear ordered non-Archimedean field $\mathcal{K}.$
%================================== definition =======================================
\begin{df} \cite{clp04}
 Let $\mathcal{X}$ and $\mathcal{Y}$ be non-Archimedean $n$-normed spaces and $f : \mathcal{X} \rightarrow \mathcal{Y}$ a mapping. We call $f$ an {\it
$n$-isometry} if
$$\|x_1 - x_0, \cdots, x_n - x_0\|=\|f(x_1) - f(x_0), \cdots, f(x_n)- f(x_0)\|$$
for all $x_0, x_1, \cdots, x_n \in \mathcal{X}$.
\end{df}

%====================== def of n-collinear =========================================================
\begin{df}~\cite{clp04}
 The points $x_{0}, x_{1}, \cdots,
x_{n}$ of a non-Archimedean $n$-normed space $\mathcal{X}$ are said to be {\it n-collinear} if for every $i$,
$\{x_{j} - x_{i} \mid 0 \le j \neq i \le n \}$ is linearly
dependent.
\end{df}

The points $x_0,x_1$ and $x_2$ of a non-Archimedean $n$-normed space $\mathcal{X}$ are said to be $2$-collinear if and only if $x_2-x_0=t(x_1-x_0)$ for some element $t$ of a real field.
We denote the set of all elements of $\mathcal{K}$ whose valuations are $1$ by $\mathcal{C},$
that is, ${\mathcal{C}}=\{\gamma\in{\mathcal{K}}:|\gamma|=1\}.$

%======================================== lemma ========================================================
\begin{lm}{\label{lm31}}
Let $x_{i}$ be an element of a non-Archimedean $n$-normed space $\mathcal{X}$ for
every $i \in \{1, \cdots , n\}$ and $\gamma\in\mathcal{K}.$
Then $$
\|x_{1}, \cdots , x_{i}, \cdots ,  x_{j}, \cdots , x_{n} \| =
\|x_{1}, \cdots , x_{i}, \cdots , x_{j} + \gamma x_{i}, \cdots ,
x_{n} \|.
$$
for all $1 \leq i \ne j \leq n.$
\end{lm}
%----------------------------------------- proof -------------------------------------------------------
\begin{pf}
By the condition $(\textrm{iv})$ of Definition~\ref{df31}, we have
\begin{eqnarray*}
&&\|x_{1}, \cdots , x_{i}, \cdots , x_{j} + \gamma x_{i}, \cdots ,x_{n} \|\\
&\leq&\max\{\|x_{1}, \cdots , x_{i}, \cdots ,  x_{j}, \cdots , x_{n} \|,|\gamma|\,\|x_{1}, \cdots , x_{i}, \cdots , x_{i}, \cdots ,x_{n} \|\}\\
&=&\max\{\|x_{1}, \cdots , x_{i}, \cdots ,  x_{j}, \cdots , x_{n} \|,0\}\\
&=&\|x_{1}, \cdots , x_{i}, \cdots ,  x_{j}, \cdots , x_{n} \|.
\end{eqnarray*}
One can easily prove the converse using the similar methods.
This completes the proof.
\end{pf}

%======================================== Remark ========================================================
\begin{rk}\label{rk33}
Let $\mathcal{X,Y}$ be non-Archimedean $n$-normed spaces over a linear ordered non-Archimedean field $\mathcal{K}$
and let $f:\mathcal{X}\rightarrow\mathcal{Y}$ be an $n$-isometry.
One can show that the $n$-isometry $f$ from $\mathcal{X}$ to $\mathcal{Y}$ preserves the $2$-collinearity
using the similar method in ~\cite[Lemma 3.2]{cck09}.
\end{rk}
%==========================================================================================================

The  {\emph{midpoint}} of a segment with endpoints $x$ and $y$ in the non-Archimedean $n$-normed space $\mathcal{X}$ is defined by the point $\frac{x+y}{2}.$

Now, we prove the Mazur-Ulam theorem on non-Archimedean $n$-normed spaces.
In the first step, we prove the following lemma. And then, using the lemma, we show that an $n$-isometry $f$ from a non-Archimedean $n$-normed space $\mathcal{X}$ to a non-Archimedean $n$-normed space $\mathcal{Y}$ preserves the midpoint of a segment.
I.e., the $f$-image of the midpoint of a segment in $\mathcal{X}$ is also the midpoint of a corresponding segment in $\mathcal{Y}.$

%============================================ Lemma ========================================================
\begin{lm}\label{lm39}
Let $\mathcal{X}$ be a non-Archimedean $n$-normed space
over a linear ordered non-Archimedean field $\mathcal{K}$ with ${\mathcal{C}}=\{2^n|\,n\in{\mathbb{Z}}\}$
and let $x_0,x_1\in\mathcal{X}$ with $x_0\neq x_1.$
Then $u:=\frac{x_0+x_1}{2}$ is the unique member of $\mathcal{X}$ satisfying
\begin{eqnarray*}
\|x_0-x_1,x_0-x_2,x_0-x_3,\cdots,x_0-x_n\|\\
=\|x_0-u,x_0-x_2,x_0-x_3,\cdots,x_0-x_n\|\\
=\|x_1-u,x_1-x_2,x_1-x_3,\cdots,x_1-x_n\|
\end{eqnarray*}
for some $x_2,\cdots,x_n\in \mathcal{X}$ with $\|x_0-x_1,x_0-x_2,\cdots,x_0-x_n\|\neq0$
and $u,x_0,x_1$ are $2$-collinear.
\end{lm}
%-------------------------------------------- proof -------------------------------------------------------
\begin{pf}
Let $u:=\frac{x_0+x_1}{2}.$
From the assumption for the dimension of $\mathcal{X},$ there exist $n-1$ elements $x_2,\cdots,x_n$ in $\mathcal{X}$ such that $\|x_0-x_1,x_0-x_2,\cdots,x_0-x_n\|\neq0.$
One can easily prove that $u$ satisfies the above equations and conditions.
It suffices to show the uniqueness for $u.$
Assume that there is an another $v$ satisfying
\begin{eqnarray*}
\|x_0-x_1,x_0-x_2,x_0-x_3,\cdots,x_0-x_n\|\\
=\|x_0-v,x_0-x_2,x_0-x_3,\cdots,x_0-x_n\|\\
=\|x_1-v,x_1-x_2,x_1-x_3,\cdots,x_1-x_n\|
\end{eqnarray*}
for some elements $x_2,\cdots,x_n$ of $\mathcal{X}$ with $\|x_0-x_1,x_0-x_2,\cdots,x_0-x_n\|\neq0$
and $v,x_0,x_1$ are $2$-collinear.
Since $v,x_0,x_1$ are $2$-collinear, $v=tx_0+(1-t)x_1$ for some $t\in\mathcal{K}.$
Then we have
\begin{eqnarray*}
&&\|x_0-x_1,x_0-x_2,x_0-x_3,\cdots,x_0-x_n\|\\
&=&\|x_0-v,x_0-x_2,x_0-x_3,\cdots,x_0-x_n\|\\
&=&\|x_0-tx_0-(1-t)x_1,x_0-x_2,x_0-x_3,\cdots,x_0-x_n\|\\
&=&|1-t|\,\|x_0-x_1,x_0-x_2,x_0-x_3,\cdots,x_0-x_n\|,
\end{eqnarray*}
\begin{eqnarray*}
&&\|x_0-x_1,x_0-x_2,x_0-x_3,\cdots,x_0-x_n\|\\
&=&\|x_1-v,x_1-x_2,x_1-x_3,\cdots,x_1-x_n\|\\
&=&\|x_1-tx_0-(1-t)x_1,x_1-x_2,x_1-x_3,\cdots,x_1-x_n\|\\
&=&|t|\,\|x_0-x_1,x_1-x_2,x_1-x_3,\cdots,x_1-x_n\|\\
&=&|t|\,\|x_0-x_1,x_0-x_2,x_0-x_3,\cdots,x_0-x_n\|.
\end{eqnarray*}
Since $\|x_0-x_1,x_0-x_2,\cdots,x_0-x_n\|\neq0$, we have two equations $|1-t|=1$ and $|t|=1.$
So there are two integers $k_1,k_2$ such that $1-t=2^{k_1},\,t=2^{k_2}.$
Since $2^{k_1}+2^{k_2}=1,$ $k_i<0$ for all $i=1,2.$
Thus we may assume that $1-t=2^{-n_1},\,t=2^{-n_2}$ and $n_1\geq n_2\in\mathbb{N}$ without loss of generality.
If $n_1\gneq n_2,$ then $1=2^{-n_1}+2^{-n_2}=2^{-n_1}(1+2^{n_1-n_2}),$ that is, $2^{n_1}=1+2^{n_1-n_2}.$
This is a contradiction because the left side of the equation is a multiple of $2$ but the right side of the equation is not.
Thus $n_1=n_2=1$ and hence $v=\frac{1}{2}x_0+\frac{1}{2}x_1=u.$
\end{pf}

%=========================================== Theorem ======================================================
\begin{thm}\label{thm38}
Let $\mathcal{X},\mathcal{Y}$ be non-Archimedean $n$-normed spaces
over a linear ordered non-Archimedean field $\mathcal{K}$ with ${\mathcal{C}}=\{2^n|\,n\in{\mathbb{Z}}\}$
and $f:\mathcal{X}\rightarrow \mathcal{Y}$ an $n$-isometry.
Then the midpoint of a segment is $f$-invariant, i.e., for every $x_0,x_1\in\mathcal{X}$ with $x_0\neq x_1,$ $f(\frac{x_0+x_1}{2})$ is also the midpoint of a segment with endpoints $f(x_0)$ and $f(x_1)$ in $\mathcal{Y}.$
\end{thm}
%----------------------------------------- proof ---------------------------------------------------------
\begin{pf}
Let $x_0,x_1\in\mathcal{X}$ with $x_0\neq x_1.$
Since the dimension of $\mathcal{X}$ is greater than $n-1,$
there exist $n-1$ elements $x_2,\cdots,x_n$ of $\mathcal{X}$ satisfying $\|x_0-x_1,x_0-x_2,\cdots,x_0-x_n\|\neq0.$
Since $x_0,x_1$ and their midpoint $\frac{x_0+x_1}{2}$ are $2$-collinear in $\mathcal{X}$,
$f(x_0),f(x_1),f(\frac{x_0+x_1}{2})$ are also $2$-collinear in $\mathcal{Y}$ by Remark~\ref{rk33}.
Since $f$ is an $n$-isometry, we have the followings
\begin{eqnarray*}
&&\|f(x_0)-f(\frac{x_0+x_1}{2}),f(x_0)-f(x_2),\cdots,f(x_0)-f(x_n)\|\\
&=&\|x_0-\frac{x_0+x_1}{2},x_0-x_2,\cdots,x_0-x_n\|\\
&=&|\frac{1}{2}|\,\|x_0-x_1,x_0-x_2,\cdots,x_0-x_n\|\\
&=&\|f(x_0)-f(x_1),f(x_0)-f(x_2),\cdots,f(x_0)-f(x_n)\|\,,\\
\\
&&\|f(x_1)-f(\frac{x_0+x_1}{2}),f(x_1)-f(x_2),\cdots,f(x_1)-f(x_n)\|\\
&=&\|x_1-\frac{x_0+x_1}{2},x_1-x_2,\cdots,x_1-x_n\|\\
&=&|\frac{1}{2}|\,\|x_1-x_0,x_1-x_2,\cdots,x_1-x_n\|\\
&=&\|f(x_0)-f(x_1),f(x_0)-f(x_2),\cdots,f(x_0)-f(x_n)\|.
\end{eqnarray*}
By Lemma ~\ref{lm39}, we obtain that $f(\frac{x_0+x_1}{2})=\frac{f(x_0)+f(x_1)}{2}$ for all $x_0,x_1\in\mathcal{X}$ with $x_0\neq x_1.$
This completes the proof.
\end{pf}

%===================================== proposition =========================================================
\begin{lm}\label{lm36}
Let $\mathcal{X}$ and $\mathcal{Y}$ be non-Archimedean $n$-normed spaces and $f:\mathcal{X} \rightarrow \mathcal{Y}$ an $n$-isometry. Then the following conditions are equivalent.

$(\textrm{i})$ The $n$-isometry $f$ is additive, i.e., $f(x_0+x_1)=f(x_0)+f(x_1)$ for all $x_0,x_1\in \mathcal{X};$

$(\textrm{ii})$ the $n$-isometry $f$ preserves the midpoint of a segment in $\mathcal{X},$
i.e., $f(\frac{x_0+x_1}{2})=\frac{f(x_0)+f(x_1)}{2}$ for all $x_0,x_1\in \mathcal{X}$ with $x_0 \neq x_1;$

$(\textrm{iii})$ the $n$-isometry $f$ preserves the barycenter of a triangle in $\mathcal{X}$, i.e., $f(\frac{x_0+x_1+x_2}{3})=\frac{f(x_0)+f(x_1)+f(x_2)}{3}$ for all $x_0,x_1,x_2\in \mathcal{X}$ satisfying that $x_0,x_1,x_2$ are not $2$-collinear.
\end{lm}
%-------------------------------------- proof ------------------------------------------------------------
\begin{pf}
Firstly, the equivalence between $(\textrm{i})$ and $(\textrm{ii})$ is obvious.
Then we suffice to show that $(\textrm{ii})$ is equivalent to $(\textrm{iii}).$
Assume that the $n$-isometry $f$ preserves the barycenter of a triangle in $\mathcal{X}.$
Let $x_0,x_1$ be in $\mathcal{X}$ with $x_0 \neq x_1.$
Since the $n$-isometry $f$ preserves the $2$-collinearity, $f(x_0),f(\frac{x_0+x_1}{2}),f(x_1)$ are $2$-collinear.
So
\begin{equation}\label{eq31}
f(\frac{x_0+x_1}{2})-f(x_0)=s\Big(f(x_1)-f(x_0)\Big)
\end{equation}
for some element $s$ of a real field.
By the hypothesis for the dimension of $\mathcal{X}$, we can choose the element $x_2$ of $\mathcal{X}$ satisfying that $x_0,x_1$ and $x_2$ are not $2$-collinear.
Since $x_2,\frac{x_0+x_1+x_2}{3},\frac{x_0+x_1}{2}$ are $2$-collinear, we have that $f(x_2),f(\frac{x_0+x_1+x_2}{3}),f(\frac{x_0+x_1}{2})$ are also $2$-collinear by Remark~\ref{rk33}.
So we obtain that
\begin{equation}\label{eq32}
f(\frac{x_0+x_1+x_2}{3})-f(x_2)=t\Big(f(\frac{x_0+x_1}{2})-f(x_2)\Big)
\end{equation}
for some element $t$ of a real field.
By the equations (\ref{eq31}), (\ref{eq32}) and the barycenter preserving property for the $n$-isometry $f$, we have $$\frac{f(x_0)+f(x_1)+f(x_2)}{3}-f(x_2)=t\Big(f(x_0)+sf(x_1)-sf(x_0)-f(x_2)\Big).$$
Thus we get $$\frac{f(x_0)+f(x_1)-2f(x_2)}{3}=t(1-s)f(x_0)+tsf(x_1)-tf(x_2).$$
So we have the following equation $$\frac{2}{3}\Big(f(x_0)-f(x_2)\Big)-\frac{1}{3}\Big(f(x_0)-f(x_1)\Big)=t\Big(f(x_0)-f(x_2)\Big)-ts\Big(f(x_0)-f(x_1)\Big).$$
By a calculation, we obtain
\begin{equation}\label{eq33}
(\frac{2}{3}-t)\Big(f(x_0)-f(x_2)\Big)+(-\frac{1}{3}+ts)\Big(f(x_0)-f(x_1)\Big)=0.
\end{equation}
Since $x_0,x_1,x_2$ are not $2$-collinear, $x_0-x_1,x_0-x_2$ are linearly independent.
Since $\dim{\mathcal{X}}\geq n,$ there are $x_3,\cdots,x_n\in \mathcal{X}$ such that $\|x_0-x_1,x_0-x_2,x_0-x_3,\cdots,x_0-x_n\|\neq0.$
Since $f$ is an $n$-isometry,
\begin{eqnarray*}
&&\|f(x_0)-f(x_1),f(x_0)-f(x_2),f(x_0)-f(x_3),\cdots,f(x_0)-f(x_n)\|\\
&&=\|x_0-x_1,x_0-x_2,x_0-x_3,\cdots,x_0-x_n\|\neq0.
\end{eqnarray*}
So $f(x_0)-f(x_1)$ and $f(x_0)-f(x_2)$ are linearly independent.
Hence, from the equation (\ref{eq33}), we have $\frac{2}{3}-t=0$ and $-\frac{1}{3}+ts=0.$
I.e., we obtain $t=\frac{2}{3}, s=\frac{1}{2},$
which imply the equation $$f(\frac{x_0+x_1}{2})=\frac{f(x_0)+f(x_1)}{2}$$ for all $x_0,x_1\in \mathcal{X}$ with $x_0\neq x_1.$

Conversely, $(\textrm{ii})$ trivially implies $(\textrm{iii})$.
This completes the proof of this lemma.
\end{pf}

%==================== remark for the equivalence relation ================================================
\begin{rk}{\label{rk3-10}}
One can prove that the above lemma also holds in the case of linear $n$-normed spaces.
\end{rk}
%==================================== Theorem ==============================================================
\begin{thm}{\label{rk37}}
Let $\mathcal{X}$ and $\mathcal{Y}$ be non-Archimedean $n$-normed spaces
over a linear ordered non-Archimedean field $\mathcal{K}$ with ${\mathcal{C}}=\{2^n|\,n\in{\mathbb{Z}}\}.$
If $f:\mathcal{X}\rightarrow \mathcal{Y}$ is an $n$-isometry, then $f(x)-f(0)$ is additive.
\end{thm}
%---------------------------------------- proof ------------------------------------------------------
\begin{pf}
Let $g(x):=f(x)-f(0).$
Then it is clear that $g(0)=0$ and $g$ is also an $n$-isometry.

From Theorem~\ref{thm38}, for $x_0,x_1\in{\mathcal{X}}(x_0\neq x_1),$  we have $$g\Big(\frac{x_0+x_1}{2}\Big)=\frac{g(x_0)+g(x_1)}{2}.$$
Hence, by Lemma~\ref{lm36}, we obtain that $g$ is additive which completes the proof.
\end{pf}

\smallskip

%================================================================================================================

In the remainder of this section, under different conditions from those previously referred in Theorem~\ref{rk37},
we also prove the Mazur-Ulam theorem on a non-Archimedean $n$-normed space.
Firstly, we show that an $n$-isometry $f$ from a non-Archimedean $n$-normed space $\mathcal{X}$ to a non-Archimedean $n$-normed space $\mathcal{Y}$ preserves the barycenter of a triangle. I.e., the $f$-image of the barycenter of a triangle is also the barycenter of a corresponding triangle.
Then, using Lemma~\ref{lm36}, we also prove the Mazur-Ulam theorem(non-Archimedean $n$-normed space version) under some different conditions.

%========================================== lemma ===========================================
\begin{lm}\label{lm35}
Let $\mathcal{X}$ be a non-Archimedean $n$-normed space over a linear ordered non-Archimedean field $\mathcal{K}$
with ${\mathcal{C}}=\{3^n|\,n\in{\mathbb{Z}}\}$ and let $x_0,x_1,x_2$ be elements of $\mathcal{X}$ such that $x_0,x_1,x_2$ are not $2$-collinear.
Then $u:=\frac{x_0+x_1+x_2}{3}$ is the unique member of $\mathcal{X}$ satisfying
\begin{eqnarray*}
\|x_0-x_1,x_0-x_2,x_0-x_3,\cdots,x_0-x_n\|\\
=\|x_0-x_1,x_0-u,x_0-x_3,\cdots,x_0-x_n\|\\
=\|x_1-x_2,x_1-u,x_1-x_3,\cdots,x_1-x_n\|\\
=\|x_2-x_0,x_2-u,x_2-x_3,\cdots,x_2-x_n\|
\end{eqnarray*}
for some $x_3,\cdots,x_n\in \mathcal{X}$ with $\|x_0-x_1,x_0-x_2,x_0-x_3,\cdots,x_0-x_n\|\neq0$
and $u$ is an interior point of $\triangle_{x_0x_1x_2}.$
\end{lm}
%------------------------------------------- proof -------------------------------------------------
\begin{pf}
Let $u:=\frac{x_0+x_1+x_2}{3}.$
Thus $u$ is an interior point of $\triangle_{x_0x_1x_2}.$
Since $\dim{\mathcal{X}}>n-1,$ there are $n-2$ elements $x_3,\cdots,x_n$ of $\mathcal{X}$ such that
$\|x_0-x_1,x_0-x_2,x_0-x_3,\cdots,x_0-x_n\|\neq0.$
Applying the lemma~\ref{lm31}, we have that
\begin{eqnarray*}
&&\|x_0-x_1,x_0-u,x_0-x_3,\cdots,x_0-x_n\|\\
&=&\|x_0-x_1,x_0-\frac{x_0+x_1+x_2}{3},x_0-x_3,\cdots,x_0-x_n\|\\
&=&|\frac{1}{3} |\,\|x_0-x_1,x_0-x_1+x_0-x_2,x_0-x_3,\cdots,x_0-x_n\|\\
&=&\|x_0-x_1,x_0-x_2,x_0-x_3,\cdots,x_0-x_n\|.
\end{eqnarray*}
And we can also obtain that
\begin{eqnarray*}
&&\|x_0-x_1,x_0-x_2,x_0-x_3,\cdots,x_0-x_n\|\\
&=&\|x_1-x_2,x_1-u,x_1-x_3,\cdots,x_1-x_n\|\\
&=&\|x_2-x_0,x_2-u,x_2-x_3,\cdots,x_2-x_n\|.
\end{eqnarray*}
For the proof of uniqueness, let $v$ be an another interior point of $\triangle_{x_0x_1x_2}$ satisfying
\begin{eqnarray*}
\|x_0-x_1,x_0-x_2,x_0-x_3,\cdots,x_0-x_n\|\\
=\|x_0-x_1,x_0-v,x_0-x_3,\cdots,x_0-x_n\|\\
=\|x_1-x_2,x_1-v,x_1-x_3,\cdots,x_1-x_n\|\\
=\|x_2-x_0,x_2-v,x_2-x_3,\cdots,x_2-x_n\|
\end{eqnarray*}
with $\|x_0-x_1,x_0-x_2,x_0-x_3,\cdots,x_0-x_n\|\neq0.$
Since $v$ is an element of the set $\{t_0x_0+t_1x_1+t_2x_2 |\,t_0+t_1+t_2=1,\,t_i\in{\mathcal{K}},\,t_i>0\,\,{\text{for all}}\,\,i\},$
there are elements $s_0,s_1,s_2$ of $\mathcal{K}$ with $s_0+s_1+s_2=1,\,s_i>0$ such that $v=s_0x_0+s_1x_1+s_2x_2.$
Then we have
\begin{eqnarray*}
&&\|x_0-x_1,x_0-x_2,x_0-x_3,\cdots,x_0-x_n\|\\
&=&\|x_0-x_1,x_0-v,x_0-x_3,\cdots,x_0-x_n\|\\
&=&\|x_0-x_1,x_0-s_0x_0-s_1x_1-s_2x_2,x_0-x_3,\cdots,x_0-x_n\|\\
&=&\|x_0-x_1,(s_0-1)x_0+s_1x_1+(1-s_0-s_1)x_2,x_0-x_3,\cdots,x_0-x_n\|\\
&=&\|x_0-x_1,(s_0+s_1-1)x_0+(1-s_0-s_1)x_2,x_0-x_3,\cdots,x_0-x_n\|\\
&=&|s_0+s_1-1|\,\|x_0-x_1,x_0-x_2,x_0-x_3,\cdots,x_0-x_n\|\\
&=&|s_2|\,\|x_0-x_1,x_0-x_2,x_0-x_3,\cdots,x_0-x_n\|
\end{eqnarray*}
and hence $|s_2|=1$ since $\|x_0-x_1,x_0-x_2,x_0-x_3,\cdots,x_0-x_n\|\neq0.$
Similarly, we obtain $|s_0|=|s_1|=1.$
By the hypothesis of $\mathcal{C},$ there are integers $k_0,k_1,k_2$ such that $s_0=3^{k_1},s_1=3^{k_2},s_2=3^{k_2}.$
Since $s_0+s_1+s_2=1,$ every $k_i$ is less than $0.$
So, one may let $s_0=3^{-n_0},s_1=3^{-n_1},s_2=3^{-n_2}$ and $n_0\geq n_1\geq n_2\in\mathbb{N}.$
Assume that the one of above inequalities holds.
Then $1=s_0+s_1+s_2=3^{-n_0}(1+3^{n_0-n_1}+3^{n_0-n_2}),$ i.e.,\,$3^{n_0}=1+3^{n_0-n_1}+3^{n_0-n_2}.$
This is a contradiction, because the left side is a multiple of $3$ whereas the right side is not.
Thus $n_0=n_1=n_2.$
Consequently, $s_0=s_1=s_2=\frac{1}{3}.$
This means that $u$ is unique.
\end{pf}

%================================================ Theorem =========================================
\begin{thm}\label{thm36}
Let $\mathcal{X},\mathcal{Y}$ be non-Archimedean $n$-normed spaces over a linear ordered non-Archimedean field $\mathcal{K}$
with ${\mathcal{C}}=\{3^n|\,n\in{\mathbb{Z}}\}$ and $f:\mathcal{X}\rightarrow\mathcal{Y}$ an interior preserving $n$-isometry.
Then the barycenter of a triangle is $f$-invariant.
\end{thm}
%--------------------------------------------- proof ---------------------------------------------------
\begin{pf}
Let $x_0,x_1$ and $x_2$ be elements of $\mathcal{X}$ satisfying that $x_0,x_1$ and $x_2$ are not $2$-collinear.
It is obvious that the barycenter $\frac{x_0+x_1+x_2}{3}$ of a triangle $\triangle_{x_0x_1x_2}$ is an interior point of the triangle.
By the assumption, $f(\frac{x_0+x_1+x_2}{3})$ is also the interior point of a triangle $\triangle_{f(x_0)f(x_1)f(x_2)}.$
Since $\dim{\mathcal{X}}>n-1,$ there exist $n-2$ elements $x_3,\cdots,x_n$ in $\mathcal{X}$ such that
$\|x_0-x_1,x_0-x_2,x_0-x_3,\cdots,x_0-x_n\|$ is not zero.
Since $f$ is an $n$-isometry, we have
\begin{eqnarray*}
&&\|f(x_0)-f(x_1),f(x_0)-f\Big(\frac{x_0+x_1+x_2}{3}\Big),f(x_0)-f(x_3),\cdots,f(x_0)-f(x_n)\|\\
&=&\|x_0-x_1,x_0-\frac{x_0+x_1+x_2}{3},x_0-x_3,\cdots,x_0-x_n\|\\
&=&|\frac{1}{3}|\,\|x_0-x_1,x_0-x_1+x_0-x_2,x_0-x_3,\cdots,x_0-x_n\|\\
&=&\|x_0-x_1,x_0-x_2,x_0-x_3,\cdots,x_0-x_n\|\\
&=&\|f(x_0)-f(x_1),f(x_0)-f(x_2),f(x_0)-f(x_3),\cdots,f(x_0)-f(x_n)\|.
\end{eqnarray*}
Similarly, we obtain
\begin{eqnarray*}
&&\|f(x_1)-f(x_2),f(x_1)-f\Big(\frac{x_0+x_1+x_2}{3}\Big),f(x_1)-f(x_3),\cdots,f(x_1)-f(x_n)\|\\
&=&\|f(x_2)-f(x_1),f(x_2)-f\Big(\frac{x_0+x_1+x_2}{3}\Big),f(x_2)-f(x_3),\cdots,f(x_2)-f(x_n)\|\\
&=&\|f(x_0)-f(x_1),f(x_0)-f(x_2),f(x_0)-f(x_3),\cdots,f(x_0)-f(x_n)\|.
\end{eqnarray*}
From Lemma~\ref{lm35}, we get $$f\Big(\frac{x_0+x_1+x_2}{3}\Big)=\frac{f(x_0)+f(x_1)+f(x_2)}{3}$$
for all $x_0,x_1,x_2\in\mathcal{X}$ satisfying that $x_0,x_1,x_2$ are not $2$-collinear.
\end{pf}

%==================== theorem related to the barycenter ================================================
\begin{thm} {\label{thm3.11}}
Let $\mathcal{X}$ and $\mathcal{Y}$ be non-Archimedean $n$-normed spaces
over a linear ordered non-Archimedean field $\mathcal{K}$ with ${\mathcal{C}}=\{3^n|\,n\in{\mathbb{Z}}\}.$
If $f:{\mathcal{X}}\rightarrow{\mathcal{Y}}$ is an interior preserving $n$-isometry, then $f(x)-f(0)$ is additive.
\end{thm}
%---------------------------------------- proof -----------------------------------------------------------
\begin{pf}
Let $g(x):=f(x)-f(0).$
One can easily check that $g(0)=0$ and $g$ is also an $n$-isometry.
Using a similar method in ~\cite[Theorem 2.4]{cck},
we can easily prove that $g$ is also an interior preserving mapping.

Now, let $x_0,x_1,x_2$ be elements of ${\mathcal{X}}$ satisfying that $x_0,x_1,x_2$ are not $2$-collinear.
Since $g$ is an interior preserving $n$-isometry, by Theorem~\ref{thm36}, $$g\Big(\frac{x_0+x_1+x_2}{3}\Big)=\frac{g(x_0)+g(x_1)+g(x_2)}{3}$$
for any $x_0,x_1,x_2\in {\mathcal{X}}$ satisfying that $x_0,x_1,x_2$ are not $2$-collinear.
From the lemma \ref{lm36},
we obtain that the interior preserving $n$-isometry $g$ is additive which completes the proof.
\end{pf}

%========================================================================================================

\medskip

\medskip

%%%%%%%%%%%%%%%%%%%%%%%%%%%%%%%%%%%%%%%%%%%%%%%%%%%%%%%%%%%%%%%%%%%%%%%%%%%%%%%%%%%%%%%%%%%%%%%%%%%%%%%%%%%%%%%%%%%%%%%%%%%%%
%%%%%%%%%%%%%%%%%%%%%%%%%%%%%%%%%%%%%%%%%%%%%%%%%%%%%%%%%%%%%%%%%%%%%%%%%%%%%%%%%%%%%%%%%%%%%%%%%%%%%%%%%%%%%%%%%%%%%%%%%%%%%

{

\end{document}